\documentstyle{article}

\def\Ent{\hbox{Ent}}
\def\Bool{\hbox{Bool}}

\def\a{\alpha}
\def\m{\mu} 
\def\b{\beta} 
\def\g{\gamma}
\def\G{\Gamma}
\def\D{\triangle}

\def\d{\delta}  
\def\th{\theta} 
\def\l{\lambda}

\def\so{\underline} 

\def\f{\rightarrow}
\def\q{\forall}
\def\e{\exists}  
\def\v{\vdash}

\def\non{\neg}

\def\c{\widehat}

\def\N{\ifmmode{\rm I\mkern-3.1mu N\mkern0.5mu}\else{\rm I\kern-.16em
N\hskip0.5pt\ }\fi\relax} 

\newtheorem{theorem}{Theorem}[section]
\newtheorem{lemma}{Lemma}[section]

\newtheorem{corollary}{Corollary}[section]

\begin{document} 

\begin{center} 
{\Large\bf  Non deterministic classical logic:

the $\l \mu^{++}$-calculus}\\[0,5cm] 
{\bf Karim NOUR}\\
LAMA - Equipe de Logique,\\
Universit\'e de Savoie\\
73376 Le Bourget du Lac\\
e-mail nour@univ-savoie.fr\\[0,5cm]
\end{center}

\begin{abstract} 
 In this paper, we present an extension of $\l \mu$-calculus
 called $\l \mu^{++}$-calculus which has the following properties:
 subject reduction, strong normalization, unicity of the
 representation of data and thus confluence only on data
 types. This calculus allows also to program the parallel-or.
\end{abstract}

\section{Introduction}

There are now many type systems which are based on classical logic ;
among the best known are the system $LC$ of J.-Y. Girard [2], the $\l
\mu$-calulus of M. Parigot [6], the $\l_c$-calculus of J.-L. Krivine
[3] and the $\lambda^{Sym}$-calculus of F. Barbanera and S. Berardi
[1]. We consider here the $\l \mu$-calculus because it has very
good properties: confluence, subject reduction and strong
normalization. On the other hand, we lose in this system the unicity
of the representation of data. Indeed, there are normal closed terms,
different from Church integers, typable by integer type (they are
called classical integers). The solutions which were proposed to solve
this problem consisted in giving algorithms to find the value of
classical integers ([5],[7]).  Moreover the presentation of typed
$\l \mu$-calculus is not very natural. For example, we do not find a
closed $\lambda \mu$-term of type $\non \non A \f A$. In this paper,
we present an extension of $\l \mu$-calculus called $\l
\mu^{++}$-calculus which codes exactly the second order classical
natural deduction. The system we propose contains a non deterministic
simplification rule which allows a program to be reduced to one of its
subroutines. This rule can be seen as a complicated {\tt garbage
collector}. This calculus which we obtain has the following
properties: subject reduction, strong normalization,
unicity of the representation of data and thus confluence only on
data types. This calculus allows also to program the parallel-or.\\

{\bf Acknowledgement.} We wish to thank C. Raffalli for helpful
discussions. We do not forget the numerous corrections and suggestions
by N. Bernard.

\section{$\l \mu$-calculus}

\subsection{Pure $\l \mu$-calculus}

$\lambda \mu$-calculus has two distinct alphabets of variables: the
set of $\lambda$-variables $x,y,z,...$, and the set of $\mu$-variables
$\alpha,\beta,\gamma$,.... Terms (also called $\lambda \mu$-terms) are
defined by the following grammar:
\begin{center}
$t :=$ $x$ $\mid$ $\lambda x\,t$
$\mid$ $(t \; t)$ $\mid$ $\mu\alpha \, [\beta]t$ 
\end{center}

The reduction relation of $\lambda \mu$-calculus is induced by fives
different notions of reduction :\\

{\bf The computation rules}
\begin{eqnarray*}
(\lambda x\, u \; v) & \rightarrow & u[x := v]\;\;\;\;\;\;\;\;(c_{\lambda}) \\
(\mu \alpha \, u \; v) & \rightarrow & \mu \alpha \,u[\alpha :=^* v] \;\;\;\;\;\;\;\;  (c_{\mu})
\end{eqnarray*}

where $u[\alpha :=^* v]$ is obtained from $u$ by replacing
inductively each subterm of the form $[\alpha]w$ by
$[\alpha](w \; v)$  \\
 
{\bf The simplification rules}
\begin{eqnarray*}
[\a]\mu \b \,u &\rightarrow &u[\b := \a] \;\;\;\;\;\;\;\;(s_1) \\
\mu \a \, [\a]u &\rightarrow &u \;\;\;\;\; (*) \;\;\;\;\;\;\;\;(s_2) \\ 
\mu \a \, u &\rightarrow &\l x \,\mu \a \, u[\alpha :=^* x] \;\;\;\;\;(**) \;\;\;\;\;\;\;\;(s_3) 
\end{eqnarray*}
(*) if $\a$ has no free occurence in $u$
 
(**) if $u$ contains a subterm of the form $[\alpha]\lambda y\,w$\\ 
 
For any $\l  \mu$-terms $t,t'$, we shall write:

-- $t \f_{\mu}^n t'$ if $t'$ is obtained from $t$ by applying $n$
times these rules.  

-- $t \f_{\mu} t'$ if there is $n \in \N$ such that $t \f_{\mu}^n t'$.\\

We have the following result ([6],[9]):

\begin{theorem}
In $\lambda \mu$-calculus, the reduction $\f_{\mu}$ is confluent.
\end{theorem} 

\subsection{Typed $\l \mu$-calculus}

Proofs are written in a second order natural deduction system with
several conclusions, presented with sequents. The connectives we use
are $\perp$, $\f$ and $\q$. We denote by $A_1,A_2,...,A_n \f A$ the
formula $A_1 \f (A_2 \f (...(A_n \f A)...))$. We do not suppose that
the language has a special constant for equality. Instead, we define
the formula $a=b$ (where $a,b$ are terms) to be $\q X\,(X(a) \f X(b))$
where $X$ is a unary predicate variable.  Let $E$ be a set of
equations. We denote by $a \approx_E b$ the equivalence binary relation
such that : if $a=b$ is an equation of $E$, then $a[x_1:=t_1,...,x_n:=t_n]
\approx_E  b[x_1:=t_1,...,x_n:=t_n]$.\\

Let $t$ be a $\lambda \mu$-term, $A$ a type, $\Gamma =
x_1:A_1,...,x_n:A_n$, $\Delta = \alpha_1:B_1,...,\alpha_m:B_m$ are two
contexts and $E$ a set of equations. The notion ``$t$ is of type $A$
in $\Gamma$ and $\Delta$ with respect to $E$'' (denoted by $\G \v t:A,\D$) is defined by the
following rules:\\

(1) $\G\v x_i :A_i,\D$ $(1\leq i\leq n)$

(2) If $\G,x:A \v t:B,\D$, then $\G\v \l x\,t:A\f B,\D$
  
(3) If $\G_1\v u:A\f B,\D_1$, and $\G_2\v v:A,\D_2$, then $\G_1,\G_2\v (u \; v):B,\D_1,\D_2$

(4) If $\G\v t:A,\D$, and $x$ not free in $\G$ and $\D$, then $\G\v t:\q x\,A,\D$

(5) If $\G\v t:\q x\,A,\D$, then, for every term $a$, $\G\v t:A[x:=a],\D$   

(6) If $\G\v t:A,\D$, and $X$ is not free in  $\G$ and $\D$, then $\G\v t:\q X\,A,\D$ 

(7) If $\G\v t:\q X\,A,\D$, then, for every formula $G$, $\G\v t:A[X:=G],\D$ 

(8) If $\G\v t:A[x:=a],\D$, and $a \approx_E b$, then $\G\v t:A[x:=b],\D$ 

(9) If $\G\v t:A,\b:B,\D$, then :

-- $\G\v \m\b\, [\a]t:B,\a :A,\D$ if  $\a \neq \b$

-- $\G\v \m\a \, [\a]t:B,\D$ if $\a = \b$ \\

The typed $\l \mu$-calculus has the following properties ([6],[8]):

\begin{theorem}\label{prop}$\;$\\
1) {\bf Subject reduction:} Type is preserved during reduction.

2) {\bf Strong normalization:} Typable $\lambda \mu$-terms are strongly normalizable. 
\end{theorem}

\subsection{Representation of data types}

Each data type generated by free algebras can be defined by a second
order formula. The type of boolean is the formula $\Bool[x] = \q X \,\{
X({\bf 1}), X({\bf 0}) \f X(x) \}$ where ${\bf 0}$ and ${\bf 1}$ are
constants. The type of integers is the formula $\Ent[x]= \q X \,\{X(0),
\q y\,(X(y) \f X(sy)) \f X(x) \}$ where $0$ is a constant symbol for
zero, and $s$ is a unary function symbol for successor.\\

In the rest of this paper, we suppose that every set of equations $E$
satisfies the following properties: ${\bf 0} \not \approx_E {\bf 1}$ 
and if $n \not = m$, then $s^n(0) \not \approx_E s^m(0)$\\

We denote by $\so{id} = \l x \,x$, $\so{{\bf 1}} = \l x \l y \,x$, $\so{{\bf 0}} = \l x \l y \,y$
and, for every $n \in \N$, $\so{n} = \l x \l y \,(y^n \; x)$ (where
$(y^0 \; x) = x$ and $(y^{k+1} \; x) = (y \; (y^k \; x))$). It is
easy to see that:

\begin{lemma}\label{lemma}$\;$\\
1) $\v \so{{\bf 1}} : \Bool[{\bf 1}]$ and $\v \so{{\bf 0}} : \Bool[{\bf 0}]$. 

2) For every $n \in \N$, $\v \so{n} : \Ent[s^n (0)]$. 
\end{lemma}

The converse of (1) lemma \ref{lemma} is true.

\begin{lemma}
 If ${\bf b} \in \{ {\bf 0} , {\bf 1} \}$ and $\v t : \Bool[{\bf b}]$, then
 $t \f_{\mu} \so{{\bf b}}$.
\end{lemma}

But the converse of (2) lemma \ref{lemma} is not true. Indeed, if we take the
closed normal term $\th=\l x\l f \,\m\a\,[\a](f \; \m\b\,[\a](f \; x))$, we have $\v \th: \Ent[s(0)]$.
 
\section{$\l \mu^{++}$-calculus}

\subsection{Pure $\l \mu^{++}$-calculus}

The set of $\l \mu^{++}$-terms is given  by the following grammar:
\begin{center}
 $t :=$ $x$ $\mid$ $\a$ $\mid$ $\l x\,t$  $\mid$ $\mu \a \, t$ $\mid$ $(t \; t)$
\end{center}
where $x$ ranges over a set $V_{\l}$ of $\l$-variables and $\a$ ranges
over a set $V_{\mu}$ of $\mu$-variables disjoint from $V_{\l}$.

The reduction relation of $\l \mu^{++}$-calculus is induced by eight
notions of reduction:\\
 
{\bf The computation rules}
\begin{eqnarray*}
(\l x  \,u \; v) &\rightharpoonup &u[x := v] \;\;\;\;\;\;\;\; (C_{\l}) \\
(\mu \a \, u \; v) &\rightharpoonup &\mu \b u\,[\a :
=\l y \,(\b \; (y \; v)) ] \;\;\;\;\;\;\;\; (C_{\mu})
\end{eqnarray*}

{\bf The local simplification rules}
\begin{eqnarray*}
((\a \;  u) \; v) &\rightharpoonup  &(\a \;  u)  \;\;\;\;\;\;\;\; (S_1) \\ 
\mu \a \mu \b \, u &\rightharpoonup &\mu \a \, u[ \b :=
{\bf id} ] \;\;\;\;\;\;\;\; (S_2)  \\
(\a \;( \b \; u)) &\rightharpoonup &(\b \;  u) \;\;\;\;\;\;\;\; (S_3) \\
(\b \; \mu \a \, u) &\rightharpoonup &u[\a := \l y \,(\b
\;  y)] \;\;\;\;\;\;\;\; (S_4)
\end{eqnarray*}
 {\bf The global simplification rules}
\begin{eqnarray*}
\mu \a \, u &\rightharpoonup &\l z \,\mu \b \, u [\a := \l y \,(\b
\; (y \; z))]  \;\;\;\;\; (*) \;\;\;\;\;\;\;\; (S_5) \\ 
\mu \a \, u[y := (\a \; v)] &\rightharpoonup &v \;\;\;\;\; (**) \;\;\;\;\;\;\;\;(S_6) 
\end{eqnarray*}
(*) if $u$ contains a subterm of the form $(\a \; \l x \,v)$

(**) if $y$ is free in $u$ and $\a$ is not free in $v$\\

For any $\l  \mu^{++}$-terms $t,t'$, we shall write 

-- $t \rightharpoonup_{\mu^{++}}^n t'$ if $t'$ is obtained from $t$ by applying
$n$ times these rules.

-- $t \rightharpoonup_{\mu^{++}} t'$ if there is $n \in \N$ such that $t \rightharpoonup_{\mu^{++}}^n t'$. \\

Let us claim first that $\l \mu^{++}$-calculus is not
confluent. Indeed, if we take $u = \l x \,\mu \a \, ((x \; (\a \; \so{{\bf
0}})) \; (\a \; \so{{\bf 1}}))$, we have (using rule $S_6$) $u
\rightharpoonup_{\mu^{++}} \l x \,\so{{\bf 0}}$ and $u
\rightharpoonup_{\mu^{++}} \l x \,\so{{\bf 1}}$.  The non confluence of
$\l \mu^{++}$-calculus does not come only from rule $S_6$. Indeed, if
we take $v = \mu \a \,((\a \; \mu \b \, \b) \so{{\bf 0}})$, we have $v
\rightharpoonup_{\mu^{++}} \mu \a \l y \,(\a \; y)$ and $v
\rightharpoonup_{\mu^{++}} \so{{\bf 0}}$. \\

The rules which are really new compared to $\l \mu$-calculus are $S_1$ and
$S_6$. The rule $S_1$ means that the $\mu$-variables are applied to
more than one term. We will see that typing will ensure this
condition. The rule $S_6$ means that if $\mu \a \,t$ has a subterm $(\a
\; v)$ where $v$ does not contain free variables which are bounded in
$\mu \a \, t$, then we can return $v$ as result. This results in the
possibility of making a parallel computation. It is clear that this
rule is very difficult to implement. But for the examples and the
properties we will present, the condition ``not active binders between
$\mu \a$ and $\a$'' will be enough. Let us explain how we can
implement the weak version of this rule. We suppose that the syntax of
the terms has two $\l$-abstractions: $\l$ and $\l'$ and two
$\mu$-abstractions: $\mu$ and $\mu'$. We write $\l'x \,u$ and $\mu' \a\,
u$ only if the variables $x$ and $\a$ do not appear in $u$. We suppose
also that for each $\mu$-variable $\a$ we have a special symbol
$\xi_{\a}$. We can thus simulate the weak version of rule $S_6$ by the
following non deterministic rules:
\begin{eqnarray*}
\mu \a \, u &\rightharpoonup &(\xi_{\a} \; u) \\
(\xi_{\a} \; \l'x \,u) &\rightharpoonup &(\xi_{\a} \; u)  \\
(\xi_{\a} \; \mu'\b \,u) &\rightharpoonup &(\xi_{\a} \; u)  \\
((\xi_{\a} \; (\a \; v)) &\rightharpoonup &v \\
((\xi_{\a} \; (u \; v)) &\rightharpoonup &(\xi_{\a} \; u) \;\;\;\;\; (*) \\
((\xi_{\a} \; (u \; v)) &\rightharpoonup &(\xi_{\a} \; v) \;\;\;\;\; (*) 
\end{eqnarray*}
(*) $u \not = \a$\\
A result of a computation is a term which does not contain symbols
$\xi_{\a}$. \\

We will see that with the exception of rule $S_6$ the $\l
\mu^{++}$-calculus is not different from $\l \mu$-calculus. We will
establish codings which make it possible to translate each one in to
the other.

\subsection{Relation between $\l \mu$- calculus and  $\l \mu^{++}$-
calculus}

We add to $\l \mu$-calculus the equivalent version of rule $S_6$:
\begin{center}
$\mu \a \, [\b] u[y := [\a]v]$ $\f'$ $v$ 
\end{center}
if $y$ is free in $u$ and $\a$ is not free in $v$. 

We denote by $\l \mu^+$-calculus this new calculus. \\

For any $\l  \mu$-terms $t,t'$, we shall write :

-- $t \f_{\mu^{+}}^n t'$ if $t'$ is obtained from $t$ by applying $n$ times these rules. 

-- $t \f_{\mu^{+}} t'$ if there is $n \in \N$ such that $t \f_{\mu^{+}}^n t'$.\\

For each $\l \mu$-term $t$ we define a $\l \mu^{++}$-term $t^*$
in the following way:
\begin{eqnarray*}
x^* &= &x \\
\{\l x \,t\}^* &= &\l x \,t^*  \\
\{(u \; v)\}^* &= &(u^* \; v^*)  \\
\{\mu \a \, [\b]t\}^* &= &\mu \a \, (\b \; t^*)
\end{eqnarray*}

We have the following result:

\begin{theorem}
Let $u,v$ be $\l \mu$-terms. If $u \f_{\mu^{+}}^n v$, then there is $m
\geq n$ such that $u^* \rightharpoonup_{\mu^{++}}^{m} v^*$.
\end{theorem}

{\bf Proof} Easy.\hfill $\Box$ \\

The converse of this coding is much more difficult to establish
 because it is necessary to include the reductions of administrative
 redexes. We first modify slightly the syntax of the $\lambda
 \mu^{++}$-calculus. We suppose that we have a particular
 $\mu$-constant $\d$ (i.e. $\mu \d \, u$ is not a term) and two other
 $\l$-abstractions: $\l^1$ and $\l^2$. The only terms build with these
 abstractions are: $\l^1 x u$ where $u$ contains only one occurence of
 $x$ and $\l^2 x x$. For the rule $C_{\mu}$, $\l$, $\l^1$ and $\l^2$
 behave in the same way. We write rules $C_{\mu}$, $S_2$, $S_4$ and
 $S_5$ in the following way:
\begin{eqnarray*}
(\mu \a \, u \; v) &\rightharpoonup &\mu \b \, u[\a : =\l^1 y\, (\b \; (y \;
v))] \;\;\;\;\;\;\;\; (C_{\mu}) \\
\mu \a \mu \b \, u &\rightharpoonup &\mu \a \, u[ \b := \l^2 x \,x ]
\;\;\;\;\;\;\;\; (S_2) \\
(\b \; \mu \a \, u) &\rightharpoonup &u[\a := \l^1 y \,(\b \;  y)]
\;\;\;\;\;\;\;\; (S_4) \\
\mu \a \, u &\rightharpoonup &\l z \mu \b \, u [\a := \l^1 y \,(\b \; (y \;
z))] \;\;\;\;\;\;\;\; (S_5)
\end{eqnarray*}

It is clear that the new $\l \mu^{++}$-calculus is stable by reductions.\\

For each $\l \mu^{++}$-term $t$ we define a $\l \mu$-term $t^{\circ}$
in the following way :
\begin{eqnarray*}
x^{\circ} &= &x  \\
\a^{\circ} &= &\l x \mu \g \, [\a] x \;\;\;\;\;(*)  \\
\{\l x \,t\}^{\circ} &= &\l x \,t^{\circ} \\
\{\l^1 x \,t\}^{\circ} &= &\l x \,t^{\circ}\\
\{\l^2 x \,x\}^{\circ} &= &\l x \mu \g \, [\d] x \\
\{\mu \a \, t\}^{\circ} &= &\mu \a \, [\d] t^{\circ} \\
\{(\l^1 x \,u \; v)\}^{\circ} &= &u^{\circ}[x := \; v^{\circ}] \\
\{(\l^2 x \,x \; v)\}^{\circ} &= &\mu \g \, [\d] v^{\circ} \;\;\;\;\;(**) \\
\{(u \; v)\}^{\circ} &= &(u^{\circ} \; v^{\circ}) \;\;\;\;\;(***) 
\end{eqnarray*}
(*) $\g \not = \a$

(**) $\g$ is not free in $v^{\circ}$

(***) $u \not = \l^i x \,w$ $i \in \{1,2\}$\\

We have the following result:

\begin{theorem}
Let $u,v$ be $\l \mu^{++}$-terms. If $u \rightharpoonup_{\mu^{++}}^n
v$, then there is $m \geq n$ and a $\l \mu$-term $w$ such that
$u^{\circ} \f_{\mu^{+}}^{m} w $ and $v^{\circ} \f_{\mu^{+}} w$.
\end{theorem}

{\bf Proof} We use the confluence of $\lambda \mu$-calculus and the following lemma:

\begin{lemma}
Let $u,v$ be $\l \mu^{++}$-terms. 

1) $\{u[x:=v]\}^{\circ} \f_{\mu^{+}} u^{\circ}[x:=v^{\circ}]$.

2) $\{u[\alpha:=\l^1 y \,(\beta \; (y
\; v))]\}^{\circ} \f_{\mu^{+}} u^{\circ}[\alpha:=^*v^{\circ}]$.\hfill $\Box$ 
\end{lemma}

We deduce the following corollary:

\begin{corollary}\label{coro}
Let $u$ be a $\l \mu^{++}$-term. If $u^{\circ}$ is strongly
normalizable then $u$ is also strongly normalizable.
\end{corollary}

\subsection{Typed $\l \mu^{++}$-calculus}

Types are formulas of second order predicate logic constructed from
$\perp$, $\f$ and $\q$. For every formula $A$, we denote
by $\neg A$ the formula $A \f \perp$ and by $\e x \,A$ the formula $\non
\q x \,\non A$. Proofs are written in the ordinary classical natural
deduction system. \\

Let $t$ be a $\l \mu^{++}$-term, $A$ a type, $\G = x_1 : A_1 ,..., x_n
: A_n, \a_1 : \non B_1 ,..., \a_m : \non B_m$ a context, and $E$ a set
of equations. We define the notion ``$t$ is of type $A$ in $\G$ with
respect to $E$'' (denoted by $\G\v' t:A$) by means of the following rules\\
 
(1) $\G\v' x_i:A_i$ $(1\leq i\leq n)$ and $\G\v'
\a_j:\non B_j$ $(1\leq j \leq m)$. 

(2) If $\G,x:A \v' u:B$, then $\G\v' \l x\,u:A \f B$. 

(3) If $\G_1\v' u:A \f B$, and $\G_2\v' v:A$, then
$\G_1,\G_2\v' (u \; v):B$. 

(4) If $\G\v' u:A$, and $x$ is not free in $\G$, then $\G\v' u:\q
x \,A$. 

(5) If $\G\v' u:\q x\,A$, then, for every term $a$, $\G\v'
u:A[x:=a]$. 

(6) If $\G\v' u:A$, and $X$ is not free in $\G$, then $\G\v' u:\q
X\,A$. 

(7) If $\G\v' u:\q X\,A$, then, for every formulas $G$, $\G\v'
u:A[X:=G]$.  

(8) If $\G\v' u:A[x:=a]$, and $a \approx_E b$, then $\G\v' u:A[x:=b]$. 

(9) If $\G,\a:\non B \v' u:\perp$, then $\G\v' \mu \a \, u:B$. \\

Consequently, we can give more explanations for rule $S_6$. It
means that ``in a proof of a formula we cannot have a subproof of the
same formula''. The terms $\mu \a \, u[y := (\a \; v)]$ and $v$ has the
same type, then the rule $S_6$ authorizes a program to be reduced to
one of its subroutines which has the same behaviour.\\

If $\D = \a_1 : B_1 ,..., \a_m : B_m$, then we denode by $\neg \D =
\a_1 : \non B_1 ,..., \a_m : \non B_m$.

If $\G = x_1 : A_1 ,..., x_n : A_n, \a_1 : \non B_1 ,..., \a_m : \non
B_m$, then we denote by $\G_{\l} =  x_1 : A_1 ,..., x_n : A_n$ and
$\G_{\mu} =  \a_1 : B_1 ,..., \a_m : B_m$.\\

We have the following results:

\begin{theorem}\label{Reltype}$\;$

1) If $\G \v t:A , \D$, then $\G , \neg \D \v' t^*:A$.

2) If $\G \v' t:A$, then  $\G_{\l} \v t^{\circ}:A , \G_{\mu} , \d : \perp$
\end{theorem}

{\bf Proof} By induction on typing. \hfill $\Box$ 

\section{Theoretical properties of $\l \mu^{++}$-calculus} 

\begin{theorem}[Subject reduction]$\;$\\
If $\G \v' u : A$ and $u \rightharpoonup v$, then $\G \v' v : A$. 
\end{theorem}

{\bf Proof} It suffices to verify that the reduction rules are well
typed. \hfill $\Box$

\begin{theorem}[Strong normalization]\label{SN}$\;$\\
If $\G \v' u : A$, then $u$ is strongly normalizable.  
\end{theorem}

{\bf Proof} According to the theorem \ref{Reltype} and the corollary
\ref{coro}, it is enough to show that the $\l \mu^+$-calculus is
strongly normalizable. It is a direct consequence of the theorem
\ref{prop} and the following lemma:

\begin{lemma}
Let $u,v,w$ be $\l \mu$-terms. If $u$ $\f'$ $v$ $\f_{\mu}^n$ $w$ then
there is $m \geq n$ and a $\l \mu$-term $v'$ such that $u$
$\f_{\mu}^m$ $v'$ $\f'$ $w$.\hfill $\Box$ 
\end{lemma}

Let $t$ be a $\l\m^{++}$-term and ${\cal V}_t$ a set of normal
$\l\m^{++}$-terms. We write $t \f_{\mu^{++}} {\cal V}_t$ iff:

-- for all $u \in {\cal V}_t$, $t \rightharpoonup_{\mu^{++}} u$.

-- If $t \rightharpoonup_{\mu^{++}} u$ and $u$ is normal, then $u \in {\cal V}_t$.

Intuitively ${\cal V}_t$ is the set of values of $t$. 

\begin{theorem}[Unicity of representation of integers]$\;$\\
If $n \in \N$ and $\v' t :\Ent[s^n (0)]$, then $t \f_{\mu^{++}} \{\so{n}\}$. 
\end{theorem}

{\bf Proof} Let $t$ be a closed normal term such that $\v' t
:\Ent[s^n(0)]$. Since we cannot use rules $S_4$ and $S_5$, we prove
that $t = \l x \l f \,u$ and $x : X(0) , f : \q y \,(X(y) \f X(s(y))) \v'
u : X(s^n(0))$. The term $u$ does not contain $\mu$-variables. Indeed,
if not, we consider a subterm $(\a \; v)$ of $u$ such that $v$ does not
contain $\mu$-variables. It is easy to see that $v$ is of the form
$(f^m \; x)$, thus $u$ is not normal (we can apply rule
$S_6$). Therefore $u = (f^n \; x)$ and $t = \so{n}$. \hfill $\Box$

\section{Some programs in $\l \mu^{++}$-calculus}

\subsection{Classical programs}

Let ${\bf \cal I} = \l x \mu \a  \, x$, ${\bf \cal C} = \l x \mu \a \,
(x
\; \a)$ and ${\bf \cal P} = \l x \mu \a \, (\a \; (x \; \a))$. It is easy
to check that:

\begin{theorem}$\;$

1) $\v' {\bf \cal I} : \q X \,\{ \perp \f X \}$, and, for every $t,t_1,...,t_n$, 
$({\bf \cal I} \; t \; t_1 ... t_n) \rightharpoonup_{\mu^{++}} \mu \a
\, t$.

2) $\v' {\bf \cal C} : \q X \,\{\non \non X \f X \}$, and, for every
$t,t_1,...,t_n$, $({\bf \cal C} \; t \; t_1 ... t_n) \rightharpoonup_{\mu^{++}}$ 

$\mu \a \, (t \; \l y \, (\a \; (y \;t_1 ... t_n )))$.

3) $\v' {\bf \cal P} : \q X \,\{(\non X \f X) \f X \}$, and, for every
$t,t_1,...,t_n$, $({\bf \cal P} \; t \; t_1 ... t_n) \rightharpoonup_{\mu^{++}}$ 

$\mu \a \, (\a \; (t \; \l y \, (\a \; (y \;t_1 ... t_n ))) \;t_1 ... t_n )$.
\end{theorem} 

Let us note that the $\l \mu^{++}$-term ${\bf \cal I}$ simulates the
{\tt exit} instruction of {\tt C} programming language and the $\l
\mu^{++}$-term ${\bf \cal P}$ simulates the {\tt Call/cc} instruction
of the {\tt Scheme} functional language (see [4]).

\subsection{Producers of integers}

For every $n_1,...,n_m \in \N$, we define the following finite
sequence $(U_k)_{1 \leq k \leq m}$: 

$U_k = (\a \; (x \; \l d \l y \, (y \; \so{n_k}) \; {\bf id} \; ({\cal I} \;
U_{k-1})))$ $(2 \leq k \leq m)$ \\
and $U_1 = (\a \; (x \; \l d \l y \, (y \; \so{n_1}) \; {\bf id} \; \a))$.

Let $P_{n_1,...,n_m} = \l x \mu \a \, U_m$. We have:

\begin{theorem}
$\v' P_{n_1,...,n_m} : \q x \,\{ \Ent[x] \f \e y \, \Ent[y] \}$, and
$(P_{n_1,...,n_m} \; \so{0}) \f_{\mu^{++}} \{ \l y \,(y \; \so{n_i})$
; $1 \leq i \leq m\}$.
\end{theorem}

{\bf Proof} For the typing, it suffices to prove that $x : \Ent[x] , \a
: \non \e y \, \Ent[y] \v' \l d \l y \, (y \; \so{n_k}) : \non \e y \,
\Ent[y] \f
\e y \, \Ent[y]$ $(1 \leq k \leq m)$ and thus $x : \Ent[x] , \a : \non \e y\,
\Ent[y] \v' U_k : \perp$ $(1 \leq k \leq m)$. 

We define the following  finite
sequence $(V_k)_{1 \leq k \leq m}$: 

$V_k = (\a \;( \l d \l y \, (y \; \so{n_k}) \;  ({\cal I} \; V_{k-1})))$
$(2 \leq k\leq m)$ and $V_1 = (\a \;  \l y \, (y \; \so{n_1}))$.

We have $(P_{n_1,...,n_m} \; \so{0}) \rightharpoonup_{\mu^{++}} \l x
\mu \a \, V_m \rightharpoonup_{\mu^{++}} \l y \, (y \; \so{n_i})$ 
$(1 \leq i \leq m)$. \hfill $\Box$ \\

Let $P_{\N} = (Y \; F)$ where 

$F = \l x \l y \mu \a \, (\a
\; (y \; \l d \, (x \; (\so{s} \; y)) \; {\bf id}  \; ({\cal I} \; (\a \; (y
\; \l d \l z \, (z \; y) \; {\bf id}  \; \a)))))$,  $Y$ is  the Turing
fixed point and $\so{s}$ a $\l \mu^{++}$-term for successor on Church
integers. It is easy to check that:

\begin{theorem}
$(P_{\N} \; \so{0}) \f_{\mu^{++}} \{ \l y \, (y \;
\so{m})$ ; $m \in \N\}$.
\end{theorem}

We can check that  $\v' F :  \q x \,\{ \Ent[x] \f \e y \, \Ent[y] \} \f  \q
x \,\{ \Ent[x] \f \e y \,\Ent[y] \}$. Therefore, if we add to the typed system
the following rule: 
\begin{center}
If $\G \v' F : A \f A$, then $\G\v' (Y \; F) : A$
\end{center}
we obtain $\v' P_{\N} : \q x \,\{ \Ent[x] \f \e y \,\Ent[y] \}$.

It is clear that, with this rule, we lose the strong normalization
property. But we possibly can put restrictions on this rule to have
weak normalization.\\

We can deduce the following corollary:

\begin{corollary}
Let ${\cal R} \subseteq \N$ be a recursively enumerable set. There is
a closed normal $\l \mu^{++}$-term $P_{\cal R}$ such that $(P_{\cal R}
\; \so{0}) \f_{\mu^{++}} \{ \so{m}$ ; $m \in {\cal R}\}$.
\end{corollary}

\subsection{Parallel-or}

Let ${\cal TB} = \{ b $ ; $b \f_{\mu^{++}} \{\so{{\bf 0}}\}$ or $b
\f_{\mu^{++}} \{\so{{\bf 1}}\} \}$ the set of true booleans.\\

A closed normal $\l \mu^{++}$-term $b$ is said to be a false
boolean iff : 

$b$ $\not \rightharpoonup_{\mu^{++}}$ $\l x \, u$ 

or

$b$ $\rightharpoonup_{\mu^{++}}$ $\l x \, u$ where $u$ $\not
\rightharpoonup_{\mu^{++}}$ $\l y \, v$ and $u$ $\not \rightharpoonup_{\mu^{++}}$
$(x \; v_1...v_n)$

or 

$b$ $\rightharpoonup_{\mu^{++}}$ $\l x \l y \, u$ where $u$ $\not
\rightharpoonup_{\mu^{++}}$ $\l y \, v$, $u$ $\not \rightharpoonup_{\mu^{++}}$
$(x \; w_1...w_n)$ and $u$ $\not \rightharpoonup_{\mu^{++}}$ $(y \;
w_1...w_n)$. \\

We denote ${\cal FB}$ the set of false booleans. Intuitively a false boolean
is thus a term which can give the first informations on a true boolean
before looping.\\

Let ${\cal B} = {\cal TB} \cup {\cal FB}$ the set of booleans.\\

We said that a closed normal $\l \m^{++}$-term $T$ is a parallel-or iff
for all $b_1,b_2 \in {\cal B}$:

$(T \; b_1 \; b_2) \f_{\mu^{++}} \{ \so{{\bf 0}} , \so{{\bf 1}}\}$ ;

$(T \; b_1 \; b_2) \rightharpoonup_{\mu^{++}} \so{{\bf 1}}$ iff $b_1 \f_{\mu^{++}} \so{{\bf 1}}$ or $b_2 \f_{\mu^{++}} \so{{\bf 1}}$ ;

$(T \; b_1 \; b_2) \rightharpoonup_{\mu^{++}} \so{{\bf 0}}$ iff $b_1 \f_{\mu^{++}} \so{{\bf
0}}$ and $b_2 \f_{\mu^{++}}  \so{{\bf 0}}$.\\

Let $or$ be a binary function defined by the following set of equations
:
 
\begin{center}
$or({\bf 1},x) = {\bf 1}$ $\;\;\;\;\;\;\;\;\;$
$or({\bf 0},x) = x$ $\;\;\;\;\;\;\;\;\;$ 
$or(x,{\bf 1}) = {\bf 1}$ $\;\;\;\;\;\;\;\;\;$
$or(x,{\bf 0}) = x$ 
\end{center} 

Let $\bigvee = \l x \l y \mu \a \, (\a \; (x \; {\bf \c{1}} \;(y \; {\bf
\c{1}} \; {\bf \c{0}}) \; ({\cal I} \; (\a \; (y \; {\bf \c{1}} \; (x
\;{\bf \c{1}} \; {\bf \c{0}}) \; \a)))))$ where ${\bf \c{1}} = \l p\,
\so{{\bf 1}}$ and ${\bf \c{0}} = \l p \,\so{{\bf 0}}$.

\begin{theorem}
$\v' \bigvee : \q x \q y \,\{ \Bool[x] , \Bool[y] \f \Bool[or(x,y)] \}$ and 
${\bigvee}$ is a parallel-or.  
\end{theorem}

{\bf Proof} Let $B[x] = \non \Bool[x] \f \Bool[x]$. 

$x : \Bool[x] \v' x : B[{\bf 1}], B[{\bf 0}] \f  B[x]$, then
$x : \Bool[x] \v' (x \;{\bf \c{1}} \; {\bf \c{0}}) :  B[x]$.

In the same way we prove that $y : \Bool[y] \v' (y \;{\bf \c{1}} \; {\bf
\c{0}}) : B[y]$.

$y : \Bool[y] \v' y : B[{\bf 1}] , B[x] \f B[or(x,y)]$, then

$x : \Bool[x] , y : \Bool[y] \v' (y \; {\bf \c{1}} \; (x \;{\bf \c{1}} \; {\bf
\c{0}})) :  \Bool[or(x,y)]$, therefore 

$\a : \non \Bool[or(x,y)], x : \Bool[x] , y : \Bool[y] \v' (\a \; (y \; {\bf \c{1}}
\; (x \;{\bf \c{1}} \; {\bf \c{0}}) \; \a))) : \perp$ and

$\a : \non \Bool[or(x,y)], x : \Bool[x] , y : \Bool[y] \v' ({\cal I} \; (\a \; (y \; {\bf \c{1}}
\; (x \;{\bf \c{1}} \; {\bf \c{0}}) \; \a)))) : \non \Bool[or(x,y)]$.

$x : \Bool[x] \v' x : B[{\bf 1}], B[y] \f B[or(x,y)]$, then

$x : \Bool[x] , y : \Bool[y] \v' (x \; {\bf \c{1}} \; (y \;{\bf \c{1}} \; {\bf
\c{0}})) :  B[or(x,y)]$, therefore 

$\a : \non \Bool[or(x,y)], x : \Bool[x] , y : \Bool[y] \v' (x \; {\bf \c{1}} \; (y \;
{\bf \c{1}} \; {\bf \c{0}}) \; ({\cal I} \; (\a \; (y \; {\bf \c{1}}
\; (x \;{\bf \c{1}} \; {\bf \c{0}}) \; \a))) :  \Bool[or(x,y)]$.

And finally : $\v' {\bigvee} : \q x \q y \,\{ \Bool[x] , \Bool[y] \f
\Bool[or(x,y)] \}$.\\

We will make three examples of reductions. Let $b_1,b_2,b_3 \in {\cal
B}$ such that $b_1 \f_{\mu^{++}} \{ \so{{\bf 0}} \}$, $b_2
\f_{\mu^{++}} \{ \so{{\bf 1}} \}$ and $b_3$
$\rightharpoonup_{\mu^{++}}$ $\l x \l y \,u$ where $u$ $\not
\rightharpoonup_{\mu^{++}}$ $\l y \,v$, $u$ $\not
\rightharpoonup_{\mu^{++}}$ $(x \; w_1...w_n)$ and $u$ $\not
\rightharpoonup_{\mu^{++}}$ $(y \; w_1...w_n)$. We will reduce
$(\bigvee \; b_1 \; b_3)$, $(\bigvee \; b_2 \; b_3)$, and $(\bigvee \;
b_3 \; b_2)$.\\

The reductions of $R_1 =(b_3 \; {\bf \c{1}} \; {\bf \c{0}})$ and $R_2
=(b_3 \; {\bf \c{1}} \; (b_i \;{\bf \c{1}} \; {\bf \c{0}}) \; \a)))$
do not terminate, and $\a$ is free in each $R$ such that $R_2
\rightharpoonup_{\mu^{++}} R$. Therefore, the only way to be compute
$(\bigvee \; b_1 \; b_3)$ and $(\bigvee \; b_2 \; b_3)$ are the
following:\\

$(\bigvee \; b_1 \; b_3)$ 

$\rightharpoonup_{\mu^{++}} \mu \a \, (\a \; ({\bf 0} \; {\bf
\c{1}} \; R'_1 \; ({\cal I} \; (\a \; R'_2 ))))$ 

$\rightharpoonup_{\mu^{++}} ...$ 

$\rightharpoonup_{\mu^{++}} \mu \a \, (\a \; ( R''_1 \; ({\cal I} \; (\a \; R''_2
))))$

$\rightharpoonup_{\mu^{++}} ...$

Then the computation does not terminate.\\

$(\bigvee \; b_2 \; b_3)$ 

$\rightharpoonup_{\mu^{++}} \mu \a \, (\a \; ({\bf 1} \; {\bf
\c{1}} \; R'_1 \; ({\cal I} \; (\a \; R'_2 ))))$ 

 $\rightharpoonup_{\mu^{++}} ...$

$\rightharpoonup_{\mu^{++}} \mu \a \, (\a \; {\bf
\c{1}} \; ({\cal I} \; (\a \; R''_2 )))$

$\rightharpoonup_{\mu^{++}} ...$

$\rightharpoonup_{\mu^{++}} \mu \a \, (\a \; {\bf
\so{1}})$ $\rightharpoonup_{\mu^{++}} {\bf \so{1}}$.\\

The reductions of $R_3 = (b_3 \; {\bf \c{1}} \; (b_2 \; {\bf \c{1}} \; {\bf
\c{0}}))$ and $R_4 = (b_3 \;  {\bf \c{1}}\;  {\bf \c{0}})$ do not terminate. 
Therefore,
the only way to compute $(\bigvee \; b_3 \; b_2)$ is the following:\\

$(\bigvee \; b_3 \; b_2)$

$\rightharpoonup_{\mu^{++}} \mu \a \, (\a \; (R'_1 \; ({\cal I} \; (\a \;
(({\bf 1} \; {\bf \c{1}}\; R'_4 )\; \a)))))$

$\rightharpoonup_{\mu^{++}} ...$

$\rightharpoonup_{\mu^{++}} \mu \a \, (\a \; (R''_1 \; ({\cal I} \; (\a \;
({\bf \c{1}} \; \a)))))$

$\rightharpoonup_{\mu^{++}} ...$

$\rightharpoonup_{\mu^{++}} \mu \a \, (\a \; (R'''_1 \; ({\cal I} \; (\a \;
{\bf \so{1}}))))$ 

$\rightharpoonup_{\mu^{++}} ...$

$\rightharpoonup_{\mu^{++}} {\bf \so{1}}$. \hfill $\Box$\\


\begin{thebibliography}{99}

\bibitem{K} F. Barbanera and S. Berardi {\em A symmetric
lambda-calculus for classical program extraction}. In M. Hagiya and
J.C. Mitchell, editors, Proceedings of theoretical aspects of computer
software, volume 789 of LNCS, pp. 495-515. Springer Verlag, 1994.

\bibitem{K} J.-Y. Girard {\em A new constructive logic: classical logic.}
Mathematical Structures in Computer Science, num 1, pp. 255-296, 1991.

\bibitem{K} J.-L. Krivine {\em Classical logic, storage operators and
2nd order lambda-calculus.} Annals of Pure and Applied Logic, num 68,
pp. 53-78, 1994.

\bibitem{K} J.-L. Krivine {\em About classical logic and imperative
programming.} Ann. of Math. and Artif. Intell., num 16, pp. 405-414, 1996. 

\bibitem{K} K. Nour {\em La valeur d'un entier classique en $\l\mu$-calcul.}
Archive for Mathematical Logic  36, pp. 461-473, 1997.

\bibitem{K} M. Parigot {\em $\lambda \mu$-calculus : an algorithm
interpretation of classical natural deduction.} 
 Lecture Notes in Artificial Intelligence, Springer Verlag, num 624,
pp. 190-201, 1992.

\bibitem{K} M. Parigot {\em Classical proofs as programs.}
Lectures Notes in Computer Science, Springer Verlag, num 713, 263-276, 1992.

\bibitem{K} M. Parigot {\em Strong normalization for second order classical natural deduction.}
Proceedings of the eighth annual IEEE symposium on logic in computer
science, pp. 39-46, 1993.

\bibitem{K} W. Py {\em Confluence en $\l \mu$-calcul.}
Th\`ese de doctorat, Universit\'e de Savoie, 1998.


\end{thebibliography}
\end{document}